\documentclass[fleqn,10pt]{wlscirep}
\usepackage[utf8]{inputenc}
\usepackage[T1]{fontenc}
\usepackage{amsmath}    
\usepackage{graphicx}  
\usepackage{amsfonts}
\usepackage{xcolor} 
\usepackage{amssymb}
\usepackage{hyperref}
\usepackage[numbers]{natbib}  %
 \newtheorem{thm}{Theorem}[section]

 \newtheorem{defn}[thm]{Definition}
 \newtheorem{rem}[thm]{Remark}
 \newtheorem{ex}{Example}
 \numberwithin{equation}{section}

\title{Optimizing Image Retrieval with an Extended {\it b}-Metric Space}

\author[1,*]{Abdelkader Belhenniche}
\author[1]{Roman Chertovskih}
\affil[1]{SYSTEC-ARISE Research Center for Systems and Technologies, Faculty of Engineering, University of Porto, Porto, Rua Dr. Roberto Frias s/n\ 4200-465, Portugal}

\affil[*]{belhenniche@fe.up.pt}


\begin{abstract}
We address the challenge of improving image retrieval in Query by Image Content (QBIC) systems, focusing on enhancing distance measures for better pattern matching. We extend the existing technique based on the non-linear elastic matching measure ${\rm NEM}_r$, satisfying the relaxed triangle inequality. We propose a new measure ${\rm NEM}_{\sigma}$, incorporating a dynamic function to adapt to varying conditions, in order to improve pattern matching accuracy. Developed within an extended $b$-metric framework, this measure ensures more flexible and accurate distance calculations. The proposed approach is given in a dynamic setting (e.g. such as moving robots), offering potential improvements in matching performance of the QBIC systems in real-world applications.
\end{abstract}
\begin{document}

\flushbottom
\maketitle

\section{Introduction}
A challenge addressed in this research is the efficient management of data storage and retrieval, especially for digital images, using the Query by Image Content (QBIC) system. A QBIC system is a database capable of handling images and searching for them based on various visual characteristics such as color, shape, and texture \cite{Flic93}. The main issue lies in establishing a mathematical framework that enhances the performance of such systems. Some of the methods, evaluating how well human perceptual differences are matched, can be found in \cite{Sca94}, while a non-linear elastic matching (NEM) distance measure was introduced in \cite{Nib95}.

In \cite{Velt2001} the authors designed a ${\rm NEM}$ distance to compare two ordered sets of contour points, $ A = \{a_1, a_2, \dots, a_m\} $ and $B = \{b_1, b_2, \dots, b_n\} $, by establishing a correspondence $f$ between the points, preserving the order. For each correspondence, a stretch $s(a_i, b_j)$ is assigned based on whether adjacent points in one set align with points in the other. The ${\rm NEM}$ distance is then computed as the minimum, in all possible correspondences, of the sum of the stretch values and the angular differences $ d(a_i, b_j) $ between the tangent angles at the corresponding points. Dynamic programming techniques were used to perform that computations. However, the considered measure does not satisfy neither the triangle inequality nor the relaxed triangle inequality, leading to inaccurate distance estimates and potential errors in pattern matching.

\begin{figure}[h]
  \centering
  \includegraphics[scale=0.45]{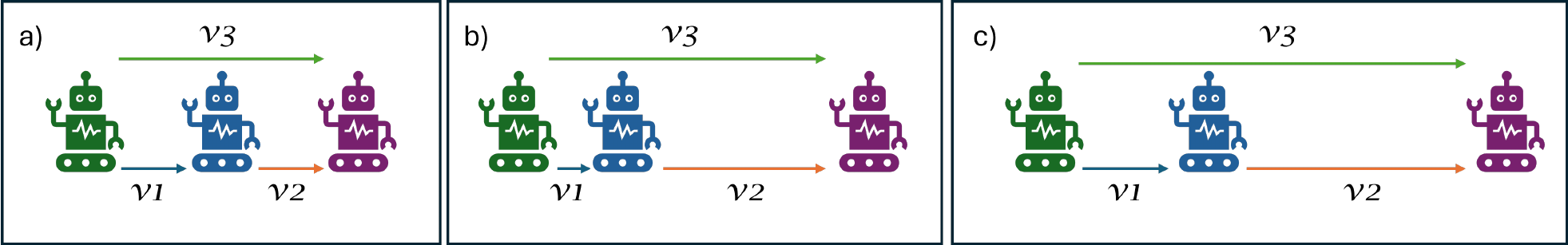}
    \caption{Pattern matching with relaxed triangle inequality. (a) represents the three robots in a stable, stationary state, while in (b) and (c) the robots are in motion.}\label{fig1}
\end{figure}

\quad To overcome this issue, an updated version of NEM, known as ${\rm NEM}_{r}$ (where $r$ is a positive constant), was developed in \cite{Nib95} and nowadays has been widely used in QBIC systems. This method involves stretching the NEM distance by a factor of $r$, ensuring that the two boundaries align. Importantly, ${\rm NEM}_{r}$ satisfies the relaxed triangle inequality for any value of $r$ and within a bounded set $S$. 


\quad In order to highlight the relevance of the relaxed triangle inequality (not related to pattern matching), we provide the following example. There are three robots (green, blue, and purple), placed in space in one straight line, see Figure \ref{fig1}(a). The robots are represented by their shapes, and the distance between them is measured as the minimal distance between their boundaries. While the distance ($v_1$) from the green robot to the blue one and the distance ($v_2$) from the blue robot to the purple are small, the distance from the green robot to the purple one is larger (because it also measures the size of the blue robot), i.e.
$$
\operatorname{NEM}(\text{Green, Blue})+
\operatorname{NEM}(\text{Blue, Purple})<
\operatorname{NEM}(\text{Green, Purple}).
$$

Note that the sign in the inequality is opposite to the sign in the classical triangle inequality, so ${\rm NEM}$ is not  even satisfy the conditions to be a metric distance. To overcome this, we may consider the ${\rm NEM}_{r}$ measure, which is a weaker form of the NEM-distance above and satisfies the following relaxed triangle inequality: 
\begin{equation}
{\rm NEM}_{r}(\text{Green, Purple}) \leq 
c\left[{\rm NEM}_{r}(\text{Green, Blue}) + 
{\rm NEM}_{r}(\text{Blue, Purple})\right]
\label{eq:bimetric}
\end{equation}
for a constant $c>1$.

\quad This approach aligns well with the theoretical foundation of $b$-metric spaces. Introduced in the seminal works \cite{Bourbaki1974} and \cite{Bakhtin1989}, and later studied in~\cite{Czerwik1993}, the $b$-metric notion relaxes the traditional triangle inequality and extends the Banach contraction mapping theorem. In this context, the use of a $b$-metric space becomes essential, as has been extensively explored in various studies. For example, in \cite{FAG98}, the authors focused on calculating the dissimilarity ratio between a given image and samples from the IBM QBIC database. They proposed methods for retrieving the most similar shape based on this ratio. Subsequently, numerous studies have emerged that address the calculation of this ratio for different forms of dissimilarity \cite{Velt2001}.

\quad In \cite{Kamran2017}, the notion of an extended $b$-metric space was introduced as a significant generalization of the standard $b$-metric space, offering a broader framework for analysis. For instance, in \cite{Karapinar2018new}, the authors derive fixed point results in an extended $b$-metric space to extend the applicability of fixed point theory to a broader range of real-life problems. They demonstrate that this framework offers a straightforward and efficient solution for the Fredholm integral equations within the context of the extended $ b $-metric space. In~\cite{panda2019new}, the authors introduce a technique related to the $F$-contraction by defining new types of contractions, namely the extended ${\mathcal{F}}_{B_{e}}$-contraction, the extended $F_{B_{e}}$-expanding contraction, and the extended generalized $F_{B_{e}}$-contraction. Based on these results, the authors proposed a simple and efficient solution for nonlinear integral equations using the fixed point technique within the framework of a $B_{e}$-metric space. Additionally, to clarify the conceptual depth of this approach, they provide illustrative examples where necessary. 


In \cite{panda2023extended}, the authors define extended suprametric spaces and establish the contraction principle using elementary properties of the greatest lower bound, departing from the traditional iteration procedure. Based on this, they derive topological results and a Stone-type theorem within the framework of suprametric spaces. Moreover, it is demonstrated that every suprametric space is metrizable. Further, they prove the existence of solutions to the Ito-Doob type stochastic integral equations using the fixed point theorem in extended suprametric spaces.

Fixed-point theory plays a crucial role in pattern matching by providing a mathematical foundation for identifying and comparing shapes based on iterative and invariant properties. In the context of Query by Image Content (QBIC) systems, fixed-point principles are used to establish correspondences between patterns or shapes, ultimately determining the closest match in a robust and consistent manner. The ability to generalize these principles across diverse conditions is essential for handling dynamic environments and complex datasets.

Frigon’s research \cite{Frigon:96} enhances this foundation by developing fixed-point results for contractive and nonexpansive mappings, by passing the need for measuring noncompactness. This simplification makes fixed-point theory more versatile and applicable in environments where traditional constraints are less practical, offering a more adaptable approach to dynamic pattern matching.

Wardowski’s work \cite{Wardowski:12} introduces a class of mappings that generalize Banach contractions, establishing a framework where fixed points exist for these mappings in complete metric spaces. This extension resolves issues related to inaccurate distance calculations, with Wardowski's $F$-contractions improving traditional methods and providing a more reliable approach to pattern matching.

Expanding on Wardowski’s results, Zaslavski \cite{Zaslavski2024} further generalizes $F$-contractions and set-valued contractions, offering more profound insights into non-expansive mappings. These developments inform the creation of the ${\rm NEM}_{\sigma}$ measure, which incorporates a dynamic function to adapt to changing conditions, ensuring more precise and flexible distance computations in real-world applications. Through these advancements, fixed-point theory provides a more comprehensive and adaptable framework for dynamic pattern matching.

In \cite{abdeljawad2019solutions}, the authors introduce a novel distance structure called the extended Branciari b-distance. This structure unifies various distance notions and yields fixed point results that encompass several existing results in the literature. As an application of their findings, they provide a solution to a fourth-order differential equation boundary value problem.

Starting from these fundamental developments, we focus on extending the capabilities of ${\rm NEM}_{r}$, which, despite its flexibility in distance computations and enhanced pattern matching accuracy, faces issues in highly complex scenarios, where additional information, such as dynamic or location-based factors, is crucial. In this article, we introduce a novel approach that extends the relaxed triangle inequality by generalizing the constant stretching penalty, $r$, to a positive function $\sigma$ defined in $X \times X$. In other words, we replace the constant $c$ in \eqref{eq:bimetric} by a dynamic function $\theta$, mapping pairs of the sequences to $[1, +\infty)$. The proposed dissimilarity measure, ${\rm NEM}_{\sigma}$, is developed within the extended $b$-metric framework. By incorporating the variable factor $\sigma(x, y)$, the measure flexibly adapts to the specific properties of the scenario under analysis. This improvement enables the computation of the distance without a recourse to a starting point on shape boundaries. Crucially, ${\rm NEM}_{\sigma}$ satisfies a generalized relaxed triangle inequality:  
\begin{equation}
   {\rm NEM}_{\sigma(x, z)}(x, z) \leq \theta(x, z) 
    \left( {\rm NEM}_{\sigma(x, y)}(x, y) + {\rm NEM}_{\sigma(y, z)}(y, z) \right), \quad \forall x, y, z \in X. 
\end{equation}\label{extendedb}
This measure is particularly effective in handling complex scenarios, such as images with varying feature distributions, dynamic objects in motion, or objects influenced by external factors such as velocity and environmental conditions.

This article is organized as follows. In Section. \ref{sec:prelim}, concepts and results for extended b-metric spaces,  ${\rm NEM}$ and  ${\rm NEM}_{r}$ measures are presented. Then, in Sect. \ref{sec:main}, as a main results of this article, we introduce the dynamic ${\rm NEM}_{\sigma}$ measure within an extended $b$-metric framework to enhance pattern matching accuracy by adapting to varying conditions for more flexible and precise distance calculations. Finally, in Sect. \ref{sec:Conclusion}, some brief conclusions and prospective research are outlined.
\section{Preliminary Concepts and Results}\label{sec:prelim}

\quad In this section we define the notation and recall the definitions, playing a key role in the derivation of further results.
\begin{defn}[\cite{Bakhtin1989, Czerwik1993}]
	Let $X$ be a non empty set, and let $c\geq 1$ be a given real number. A functional $b\colon X\times X\to [0,\infty)$ is said to be a $b$-metric if the following conditions are satisfied:
	\begin{enumerate}
		\item $b(x,y)=0$ if and only if $x=y$,
		\item $b(x,y)=b(y,x)$,
		\item $b(x,z)\leq c[b(x,y)+b(y,z)]$,
	\end{enumerate}
	for all $x,y,z\in X$. A pair $(X,b)$ is called a $b$-metric space.
\end{defn}
 Numerous fixed-point theorems in $b$-metric spaces have been developed and their applications thoroughly investigated (see, e.g. \cite{Bakhtin1989, Czerwik1993, belhenniche:22,  Kirk:01, Rao:13, Roshan:13}). 

 \begin{ex}[\cite{Roshan:13}]
	Let $(X,b)$ be a metric space, and $\rho(x,y) = (b(x,y))^{p}$, where $p\geq 1$ is a real number. Then, $(X,\rho)$ is a $b$-metric space with $\displaystyle c=2^{p-1}$.
\end{ex}
It is clear that a $b$-metric space becomes a metric space if we take $c=1$. This clearly shows that the class of $b$-metric spaces is larger than that of metric spaces.


As we already mentioned in introduction, the concept of an extended $b$-metric space is introduced in \cite{Kamran2017}, generalizing the traditional $b$-metric space and providing a more flexible framework for the analysis. In the following, we present the definition of an extended $b$-metric space:
\begin{defn}
    Let $ X $ be nonempty, and $ \theta\colon X \times X \rightarrow [1,+\infty)$. A function $ b_{\theta}\colon X \times X \rightarrow [0, +\infty) $ is an extended $b$-metric if, for all $ x, y, z \in  X $, it satisfies:
	\begin{enumerate}
		\item[1)] $ b_{\theta}(x, y)=0$ if and only if  $x=y$
		\item[2)] $ b_{\theta}(x, y)=b_{\theta}(y, x)$
		\item[3)] $ b_{\theta}(x, z) \leq \theta(x, z) [b_{\theta}(x, y)+b_{\theta}(y, z)]$
	\end{enumerate}
	The pair $ (X, b_{\theta}) $ is called extended b-metric space.\\
	Denote the open ball (the closed ball, respectively) of radius $ l >0 $ about $ x $ as the set 
	$$ B_{l}(x) =\{ x \in X\colon  b_{\theta}(x, y) < l\}, \;( B_{l}[x] =\{ x \in X\colon  b_{\theta}(x, y) \leq l\}).$$
\end{defn}

\begin{rem} If $ \theta (x, y)=c  $ for $ c \geq  1 $, then $(X,b_\theta)$ satisfies the definition of a $ b $-metric space.
\end{rem}
\begin{ex}
	Let $ X=[0, +\infty) $, and mappings $b_\theta$ and $ \theta$ with $\theta\colon X \times X \to [1, +\infty)$, defined by
	$ b_{\theta} (x, y) = (x-y)^{2}$ and $\theta(x, y)= x+y+2$.
	Then, $(X, b_{\theta})$ is an extended $b$-metric space.
\end{ex}
\begin{ex} Let $X= C([a,b], \mathbb{R})$ be the space of all continuous real valued functions defined on $ [a, b] $. Let $\displaystyle b_{\theta} (x, y)=\sup_{t \in [a, b]} \{\| x(t) - y(t) \|^{2}\}$, and $ \theta\colon X \times X \rightarrow [1, +\infty) $ defined by 
$$ \theta (x, y):= \| x(t) \| + \| y(t) \| +2, $$
then $ (X , b_\theta)$ is a complete extended $ b $-metric space.
\end{ex}

\begin{defn}
	Let $ (X, b_{\theta}) $ be an extended $ b $-metric space.
	\begin{itemize}
		\item[(i)] A sequence $\displaystyle \{ x_{n} \}_{n \in \mathbb{N}}$ in $ X $ converges to $ x \in X $ if, for every $ \varepsilon > 0 $, there exists $ N=N(\varepsilon)  \in \mathbb{N}$ such that $$ b_{\theta}(x_{n}, x) < \varepsilon $$ for all $ n \geq N $. Alternatively we may write $\displaystyle \lim_{n \rightarrow \infty} x_{n}=x $.
		\item[(ii)] A sequence $\displaystyle \{ x_{n} \}_{n \in \mathbb{N}}$ in $ X $ is Cauchy, if for every $ \epsilon > 0 $, there exists $ N=N(\varepsilon) \in \mathbb{N} $ such that $$ b_{\theta} (x_{m}, x_{n}) < \varepsilon,$$ for all $m, n \geq  N$.
	\end{itemize}
\end{defn}
\begin{defn} An extended $b$-metric space $ (X, b_{\theta})$ is complete if every Cauchy sequence in $ X $ is convergent.
\end{defn}


We now shift focus to the nonlinear elastic matching  ${\rm NEM}$  distance, a measure discussed in \cite{Velt2001}. According to this approach, given two finite sets of ordered contour points, $A = \{a_1, \dots, a_m\}$ and $B = \{b_1, \dots, b_n\}$, and a correspondence $f$ between all points in $A$ and $B$, it is required that no indices $a_1 < a_2$ exist such that $f(a_1) > f(a_2)$. The stretch $s(a_i, b_j)$ of the pair $(a_i, f(a_i) = b_j)$ is set to 1 if either $f(a_{i-1}) = b_j$ or $f(a_i) = b_{j-1}$, and 0 otherwise.

The distance ${\rm NEM}(A, B)$ is defined as the minimum over all correspondences of the sum:
$${\rm NEM}(A, B)= \min_{f} \left( \sum_{i=1}^m \sum_{j=1}^n s(a_i, b_j) + d(a_i, b_j) \right), $$
where $d(a_i, b_j)$ measures the difference between the tangent angles at $a_i$ and $b_j$. This distance can be computed efficiently using dynamic programming~\cite{Cor:94}. However, it is important to note that the NEM distance does not satisfy the triangle inequality, so it is not metric.


Turning to the ${\rm NEM_r}$ distance, \cite{FAG98} introduced the relaxed triangle inequality for this measure, as expressed in \eqref{eq:bimetric}. The authors demonstrated that the constant $c$ in the inequality is independent of the sequence lengths of $A$, $B$, and $C$. In the context of shape matching, this independence implies that $c$ does not depend on the number of sample points. When the sample points are uniformly distributed across all shapes, $c$ is given by $1 + \pi/{2r}$, while for shapes with varying sample points, $c$ becomes $1 + \pi/{r}$. Unlike weaker formulations of the relaxed triangle inequality, which allow $c$ to depend on the dimensionality of the space, such a dependency is unnecessary for the ${\rm NEM_r}$.


 Further explanation of the ${\rm NEM_r}$ distance measure is provided in \cite{FAG98}, describing it in terms of a mapping $M$, where each pair $\langle i, j \rangle$ represents an edge with certain conditions. Every integer from 1 to $m$ must appear as the first component of some edge, and similarly, every integer from 1 to $n$ must appear as the second component. The edges must not cross, meaning there should be no pairs $\langle i, j' \rangle$ and $\langle i', j \rangle$ where $i < i'$ and $j < j'$. A mapping is minimal if no smaller subset of edges can form a valid $m \times n$ mapping. Stretched edges in $M$ are those where either $\langle i-1, j \rangle$ or $\langle i, j-1 \rangle$ also exists in the mapping, and these edges are assigned a stretch cost of $r$, with non-stretched edges having a cost of 0. The distance cost between any two points is determined by the function $b(x_i, y_i)$, which measures the difference between the corresponding elements in the sequences $X$ and $Y$. The total cost of the mapping is the sum of all stretch and distance costs for the edges.

\begin{figure}[h]
  \centering
  \includegraphics[scale=0.7]{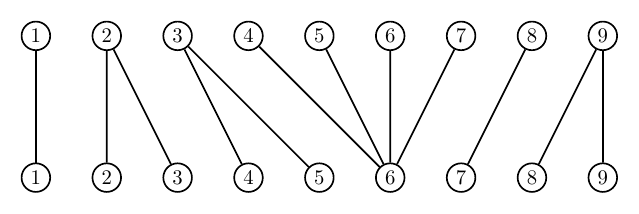}
    \caption{A minimal $(9,9)$-mapping with the stretch-cost of $6r$ (reproduced from \cite{FAG98}, Fig.~1).\label{fig2}
    }
\end{figure}

To illustrate this, we present in Figure \ref{fig2} (taken from \cite{FAG98}) the edges $\langle 2, 3 \rangle$, $\langle 3, 5 \rangle$, $\langle 5, 6 \rangle$, $\langle 6, 6 \rangle$, $\langle 7, 6 \rangle$, and $\langle 9, 9 \rangle$, which are identified as stretched edges, each incurring a stretch cost of $r$, while all other edges have a stretch cost of~$0$.

\begin{figure}[h]
  \centering
  \includegraphics[scale=0.35]{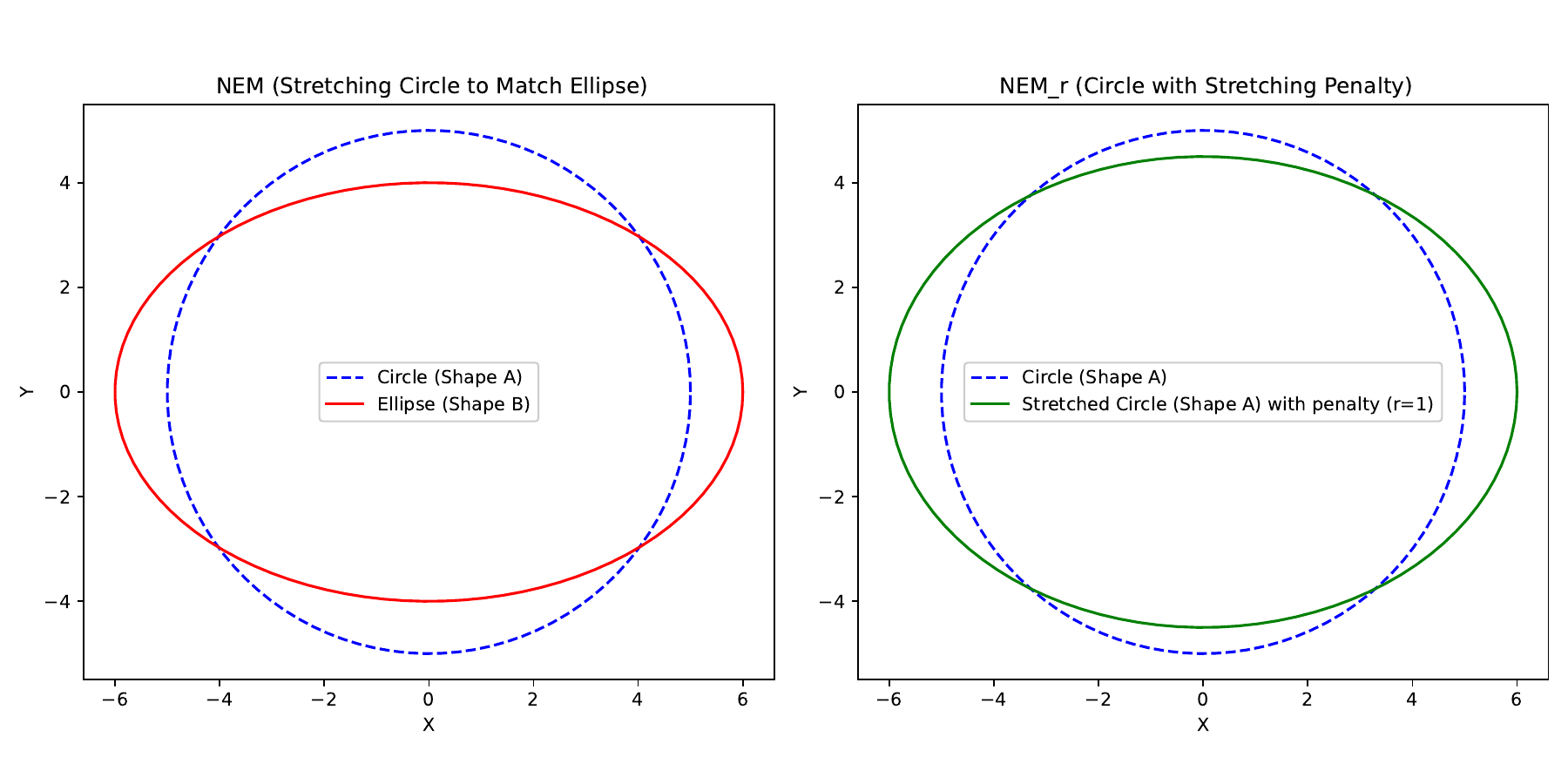}
    \caption{Comparison of ${\rm NEM}$ and ${\rm NEM}_r$: The left plot illustrates the alignment of a circle and an ellipse in ${\rm NEM}$, where the circle remains unchanged. In contrast, the right plot demonstrates ${\rm NEM}_r$, where a stretching penalty ($r = 1$) elongates the circle along the $y$-axis, modifying its shape to better match the dimensions of the ellipse. }
    \label{fig3} 
\end{figure}

To compare ${\rm NEM}$ and ${\rm NEM}_r$ effectively, we refer to Figure \ref{fig3}, showing a circle and an ellipse. It illustrates the stretching in ${\rm NEM}_r$ and demonstrates how adding a stretching penalty improves the accuracy and sensitivity of the distance measurement in ${\rm NEM}_r$.
 In the left panel,  ${\rm NEM}$ allows the circle (blue dashed line) to stretch freely to match the ellipse (red line), with no penalty for stretching. The distance between the shapes is based purely on how well their boundaries align. In the right panel, ${\rm NEM}_r$ applies a penalty for excessive stretching. The green line represents the stretched circle, and the parameter $r$ controls the amount of penalty. This ensures that extreme deformations are penalized, leading to a more realistic and meaningful distance measure between the shapes.

\begin{figure}[h!]
  \centering
  \includegraphics[scale=0.4]{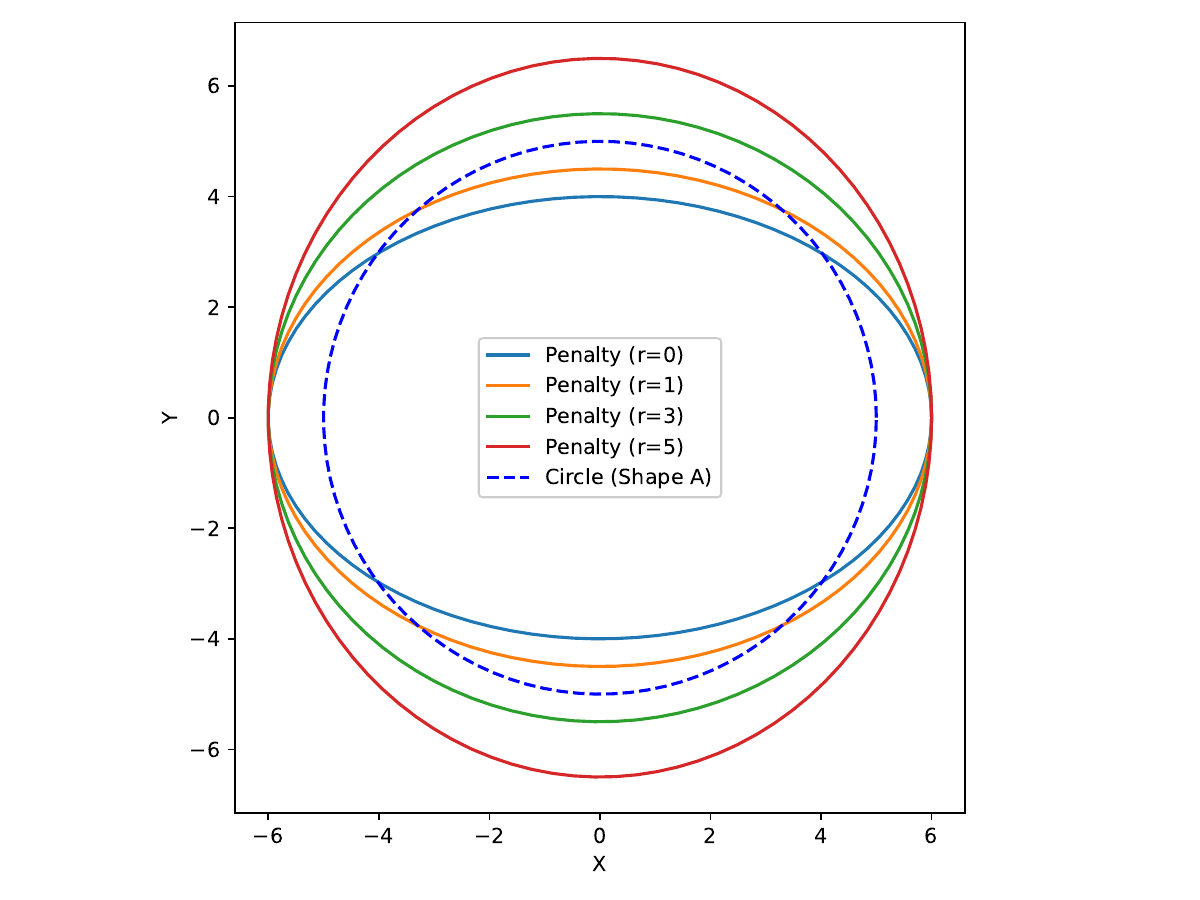}
    \caption{  Effect of parameter $r$ on stretching penalty in shape matching}\label{fig4}
\end{figure}
In figure \ref{fig4}, the parameter $r$ controls the stretching penalty during shape matching. When $r = 0$, there is no penalty, and the circle deforms minimally. As $r$ increases, the stretching penalty grows, leading to more significant deformation. A higher $r$ discourages excessive stretching, ensuring that the circle remains closer to its original shape while matching the ellipse, resulting in a more realistic and controlled distance measurement.

\section{Main results}\label{sec:main}
\quad The scenario illustrated in Figure \ref{fig1} (b) and (c) serves as a key motivation for our main results. In this case, we consider robots moving with variable velocities while remaining aligned along the same line, as depicted in Figure \ref{fig1} (a). As the robots approach and move apart, the distances between them evolve into a complex function of their initial positions and velocities. In such a dynamic setting, the relaxed triangle inequality \eqref{extendedb} continues to hold, but it is now influenced by a function $\theta(\cdot,\cdot)$, which accounts for the velocities of the objects. Furthermore, given the complexity introduced by the velocity-dependent motion, we extend the analysis by utilizing a function $\sigma(\cdot,\cdot)$ instead of a simple distance measure $r$, which provides a more accurate characterization of the stretching behavior. This framework, which blends velocity-dependent adjustments and refined distance measures, forms the foundation of our main results, which revolve around these ideas and their implications for optimal matching.

\begin{figure}[h]
  \centering
  \includegraphics[scale=0.4]{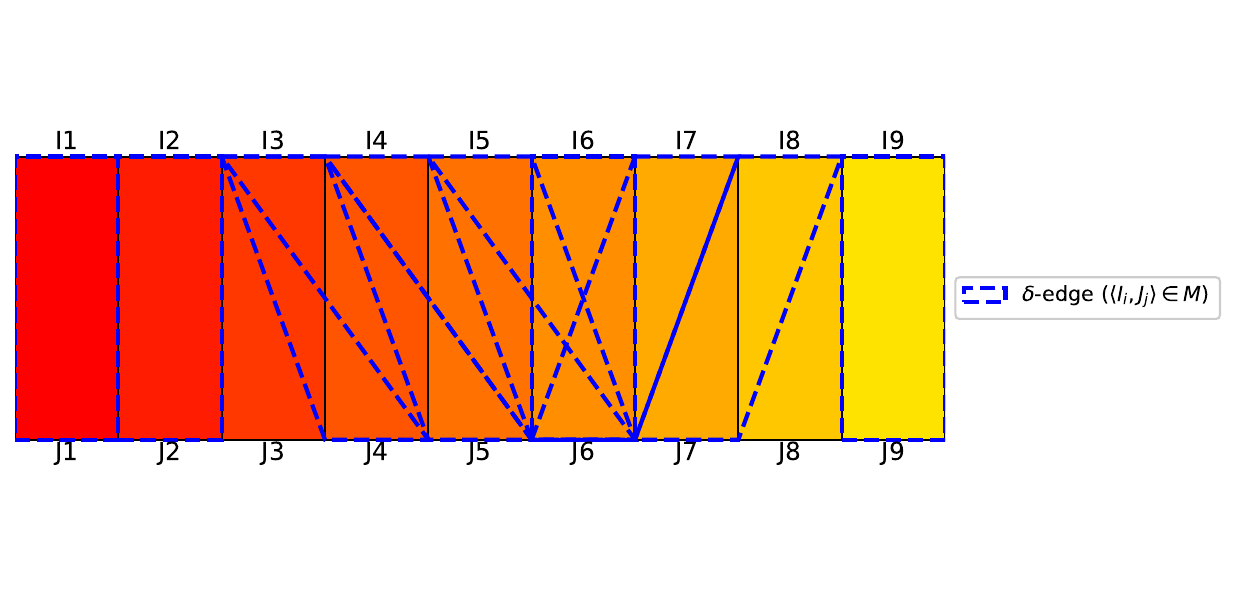}
  \caption{Minimal (9,9)-$\sigma$-mapping.}
  \label{fig5}
\end{figure}

The Figures \ref{fig2}, \ref{fig5} demonstrate how the $(9,9)$-mapping framework evolves from the rigid structure of $ \text{NEM}_r $ (Figure \ref{fig2}) to the more flexible $ \text{NEM}_\sigma $ (Figure \ref{fig5}) through the extension of the $b$-metric space. In Figure \ref{fig2}, $ \text{NEM}_r $ relies on the standard $b$-metric, using integer-based edges $\langle i, j \rangle$ with fixed stretch costs $r$ for stretched edges, resulting in a total stretch cost of $6r$. Figure \ref{fig5} extends this framework into the $ \text{NEM}_\sigma $ domain by adopting an extended $b$-metric space, where the intervals are subdivided into $\delta$-edges, and stretch costs are dynamically computed as $ \sigma(x_i, y_j) $. This extension provides greater flexibility, capturing more nuanced alignment features such as deformations and velocity adjustments. By transitioning from a discrete $b$-metric to an extended $b$-metric space, $ \text{NEM}_\sigma $ offers a more adaptable and precise model for complex sequence matching scenarios.

Building on the dynamic scenario described in the previous section, we now shift our focus to demonstrating how an extended $b$-metric space can be applied for optimal pattern matching. The goal of matching two sequences is to find the minimal combined cost, which includes both the stretch and distance costs. This concept is closely related to the idea of adjustments based on velocity, where the distances between elements evolve as a complex function of their initial positions and velocities, as highlighted earlier. To implement this approach, we first recall the essential idea of subdividing a segment, which is crucial for defining the edges in the matching process. 

Consider a compact interval $[1, n]$, and let $P = \{x_0, x_1, \ldots, x_n\}$ represent a finite set of real numbers where $1 = x_0<x_1<\ldots<x_n = n$. This subdivision of the interval generates smaller subintervals, which we define as:
$$ I_1 = [x_0, x_1], \quad I_2 = [x_1, x_2], \quad \ldots, \quad I_n = [x_{n-1}, x_n].$$
The length of each subinterval, denoted by $\delta x_i = x_i - x_{i-1}$, reflects the resolution of the segmentation. This process allows us to break down the intervals into discrete components for comparison.

Next, we consider two compact intervals, $[1, m]$ and $[1, n]$, which have both been subdivided. An $(m, n)$-$\sigma$-mapping is a set $M \subseteq [1, m] \times [1, n]$. Each pair $\langle I_i, J_j \rangle \in M$ is called a $\delta$-edge if it satisfies the following conditions:
\begin{enumerate}
    \item Every subdivision in $[1, m]$ is the first component $I_i$ of some $\delta$-edge $\langle I_i, J_j\rangle \in M$.
    \item Every subdivision in $[1, n]$ is the second component $J_j$ of some $\delta$-edge $\langle I_i, J_j \rangle \in M$.
    \item There do not exist $I_i, I_i', J_j, J_j'$ with $\delta x_i < \delta x_i'$ and $\delta y_j < \delta y_j'$ that would violate the ordering.
\end{enumerate}

\begin{figure}[h]
  \centering
  \includegraphics[scale=0.4]{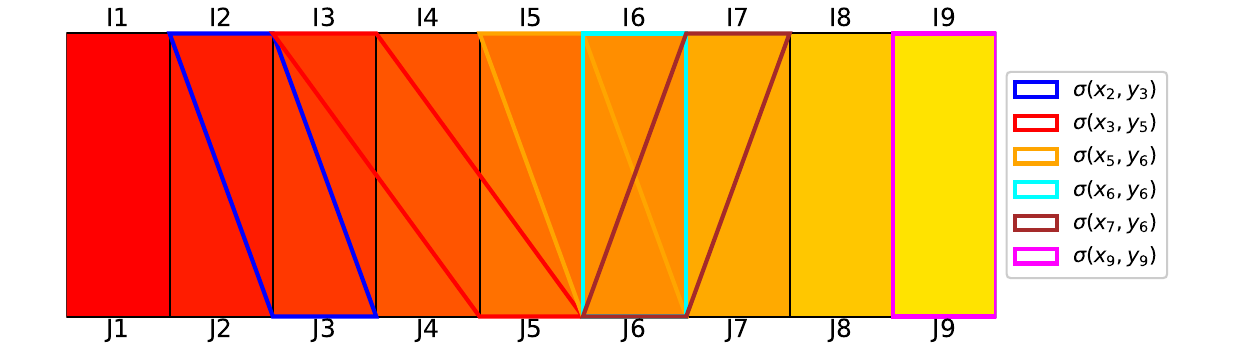}
  \caption{The Stretch Cost Representing the Minimal $(9,9)$-$\sigma$-Mapping Depicted in Figure 6:
$\sigma(x_2, y_3) + \sigma(x_3, y_5) + \sigma(x_5, y_6) + \sigma(x_6, y_6) + \sigma(x_7, y_6) + \sigma(x_9, y_9)$
}
  \label{fig6}
\end{figure}

A $\delta$-edge is a  $\sigma$-stretch-edge if both $\langle I_{i-1}, J_j \rangle$ and $\langle I_i, J_{j-1} \rangle$ are present in~$M$.

\subsection*{Cost Definitions}

For each $\delta$-edge $\langle I_i, J_j \rangle \in M$, we define two types of costs:
\begin{itemize}
    \item The stretch cost denoted by $\mathcal{S}$-cost$(\langle I_i, J_j \rangle, M)$, measures how much the pattern is stretched during matching. It is defined as:
    $$
    \mathcal{S} \text{-cost}(\langle I_i, J_j \rangle, M) =
    \begin{cases}
        \sigma(x_i, y_j), & \text{if } \langle I_i, J_j \rangle \text{ is a } \sigma\text{-stretch-edge}, \\
        0, & \text{otherwise}.
    \end{cases}
    $$
    For instance, in Figure \ref{fig5}, the minimal (9,9)-$\sigma$-mapping is depicted, highlighting specific $\delta$-edges as $\sigma$-stretch edges. The $\sigma$-stretch edges include $\langle I_2, J_3 \rangle$, $\langle I_3, J_5 \rangle$, $\langle I_5, J_6 \rangle$, $\langle I_6, J_6 \rangle$, $\langle I_7, J_6 \rangle$, and $\langle I_9, J_9 \rangle$. The total cost associated with these $\sigma$-stretch edges is given by the expression
    $$\sigma(x_2, y_3) + \sigma(x_3, y_5) + \sigma(x_5, y_6) + \sigma(x_6, y_6) + \sigma(x_7, y_6) + \sigma(x_9, y_9).$$
    All other $\delta$-edges in the mapping are associated with a stretch cost of zero.

    \item The distance cost, denoted by $\mathcal{D}$-cost$(\langle I_i, J_j \rangle)$, measures the distance between corresponding subdivisions, i.e. 
    $ \mathcal{D}$-cost$(\langle I_i, J_i \rangle) = \mathcal{B}(x_i, y_j)$, where $ \mathcal{B}(x_i, y_j) $  represents the cost associated with aligning segment $I_i $ from sequence $X$  with segment $J_j $ from sequence $ Y $. This measure captures various alignment challenges, including shape deformation, feature mismatches, and other discrepancies. By quantifying these alignment costs, $\mathcal{B}(x_i, y_j) $ provides a flexible and comprehensive way to evaluate how well different segments of sequences fit together, making it an essential component of extended $b$-metric spaces.

\end{itemize}
\section*{Total Cost of the $(m, n)$-$\sigma$-Mapping}

The total cost of the $(m, n)$-$\sigma$-mapping $M $ is the sum of two components: the stretch cost and the distance cost. Let us define these costs more precisely.

The total stretch cost is given by the double integral over the stretch function $ \sigma(x, y) $, which represents the stretching penalty when matching the elements $ x $ and $ y$. It is expressed as:

$$\mathcal{S}-cost(M,X,Y)= \int_{I} \int_{J} \sigma(x, y) \, \mathrm{d}x \, \mathrm{d}y,$$
where $\sigma(x, y) $ represents how much $x $ has to be stretched to match $y_j $.

The total distance cost is given by the double integral over the distance function $ \mathcal{B}(x, y) $, which quantifies the distance between $x$ and $y $. This cost is expressed as:

$$ \mathcal{D}-cost(M,X,Y)=\int_{I} \displaystyle\int_{J} \mathcal{B}(x, y) \, \mathrm{d}x \, \mathrm{d}y.$$

Thus, the overall total cost of the mapping \( M \), denoted as \( Cost(M, X, Y) \), is the sum of both the stretch and distance costs:

$$\text{Cost}(M, X, Y) = 
\mathcal{D}-cost(M,X,Y)\ +\ 
\mathcal{S}-cost(M,X,Y).
$$

The  extended matching cost, denoted by ${\rm NEM}_{\sigma(x, y)} $, is the minimum total cost over all possible $(m, n)$-$\sigma$-mappings \( M \). That is,

$$
{\rm NEM}_{\sigma(x, y)}(x, y) = \min_M Cost(M, X, Y).
$$

This cost represents the optimal matching cost between the sequences $X$  and $Y$, considering both the stretch and distance costs. In other words, ${\rm NEM}_{\sigma(x,y)}$  captures the minimal cost required to align the sequences $X$ and $Y$ under the given $\sigma $-mapping, taking into account all possible pairings of elements from the two sequences.

\begin{thm}
 Let $ X $ be nonempty sequence and let the function $ \mathcal{B} \colon X \times X \rightarrow [0, +\infty) $ be an extended $b$-metric with modulus  $ \alpha \colon X \times X \rightarrow [1,+\infty)$. Then, for all $ X, Y, Z$, the measure ${\rm NEM}_{\sigma} $  is extended b-metric.
\end{thm}
\textbf{Proof}:
The symmetry of the matching cost is straightforward. For any sequences $x$ and $y$, both the stretch cost $\sigma(x_i, y_j)$ and distance cost $\mathcal{B}(x_i, y_j)$ are symmetric, meaning $\mathcal{B}(x_i, y_j) = \mathcal{B}(y_j, x_i)$ and $\sigma(x_i, y_j) = \sigma(y_j, x_i)$. Consequently, the total matching cost satisfies:

$$ {\rm NEM}_{\sigma(x, y)}(x, y) = {\rm NEM}_{\sigma(y, x)}(y, x). $$

Thus, the symmetry property holds.


\quad The matching cost is non-negative because both $\sigma(x_i, y_j) \geq 0$ and $b(x_i, y_j) \geq 0$ for all $(x_i, y_j)$. Hence:

$${\rm NEM}_{\sigma(x, y)}(x, y) \geq 0.$$

For the identity of indiscernibles, if $x = y$, both $\mathcal{B}(x_i, y_j)$ and $\sigma(x_i, y_j)$ are zero, as no distance or stretching is needed, yielding:

$${\rm NEM}_{\sigma(x, y)}(x, y) = 0.$$

Conversely, if ${\rm NEM}_{\sigma(x, y)}(x, y) = 0$, both costs must be zero for all pairs $(x_i, y_j)$, which implies $x = y$. Thus, we have:

$$
{\rm NEM}_{\sigma(x, y)}(x, y) = 0 \quad \text{if and only if} \quad x = y.
$$

\medskip

Now, let us prove the relaxed triangle inequality. Starting from the definition:  
$$
{\rm NEM}_{\sigma(x, z)}(x, z) = \min_M Cost(M, X, Z),
$$ 
we express this cost as:  

\begin{equation}\label{eq: relax1}
 Cost(M, X, Z)= \mathcal{S}-cost(M,X,Z) +\mathcal{D}-cost(M,X,Z)
\end{equation}

More precisely:  
$$ Cost(M, X, Z) = \int_{I} \int_{K} (\sigma(x, z) + \mathcal{B}(x, z)) \, \mathrm{d}x \, \mathrm{d}z,$$ 
Since $\mathcal{B}$ is an extended b metric measure we get: 
\begin{equation}\label{eq: relax2}
   \mathcal{B}(x, z) \leq \alpha(x, z) [\mathbf{B}(x, y) + \mathcal{B}(y, z)]  
\end{equation}
Applying the double integralto both sides of the inequality \eqref{eq: relax2}, we get:  
$$
\int_{I} \int_{K} \mathcal{B}(x, z) \leq \Delta(x, z) \left( \int_{I} \int_{J} \mathcal{B}(x, y) \, \mathrm{d}x \, \mathrm{d}y + \int_{J} \int_{K} \mathcal{B}(y, z) \, \mathrm{d}y \, \mathrm{d}z \right),  
$$
where $$\Delta(x, z) = \int_{I} \int_{K} \alpha(x, z) \, \mathrm{d}x \, \mathrm{d}z > 1$$. Thus, we have:  
\begin{equation}\label{eq: relax3}
  \mathcal{D}-cost(M,X,Z)     \leq \Delta(x,z) (\mathcal{D}-cost(M,X,Y)  +   \mathcal{D}-cost(M,Y,Z) ),   
\end{equation}

Similarly, the strech cost satisfies:  
\begin{equation}\label{eq: relax4}
    \mathcal{S}-cost(M,X,Z) \leq   \mathcal{S}-cost(M,X,Y)+   \mathcal{S}-cost(M,Y,Z),
\end{equation}
since  
$$
\int_{I} \int_{K} \sigma(x, z) \, \mathrm{d}x \, \mathrm{d}z \leq \int_{I} \int_{J} \sigma(x, y) \, \mathrm{d}x \, \mathrm{d}y + \int_{J} \int_{K} \sigma(y, z) \, \mathrm{d}y \, \mathrm{d}z.
$$

Summing  \eqref{eq: relax3} with \eqref{eq: relax4} and using the total cost \eqref{eq: relax1}, we obtain:  
$$
Cost(M, X, Z) \leq (1 + \Delta(x, z)) [Cost(M, X, Y) + Cost(M, Y, Z)].
$$

Taking the infimum over all mappings \(M\), we conclude:  
$$
{\rm NEM}_{\sigma}(x, z) \leq \theta(x, z) \left( {\rm NEM}_{\sigma}(x, y) + {\rm NEM}_{\sigma}(y, z) \right),
$$
where $\theta(x, z) = 1 + \Delta(x, z)$. This completes the proof.

\section{Conclusion}\label{sec:Conclusion}
\quad In this article, we have introduced the ${\rm NEM}_{\sigma}$ measure as a novel technique for Query by Image Content systems, leveraging the powerful framework of extended $b$-metric spaces. By constructing a distance measure that satisfies a relaxed version of the triangle inequality, our ultimate goal is to construct a tool for more accurate and efficient image retrieval, even for non-steady objects. This new measure has the potential to increase the precision of search results,  improving the scalability and robustness of data management systems. In this paper, we have presented an analysis of some mathematical properties of the measure. However, our findings also suggest that ${\rm NEM}_{\sigma}$ can be a useful tool in the advancement of practical image retrieval technologies. As a next step of our research, we plan to validate the effectiveness of ${\rm NEM}_{\sigma}$ in some real-world applications, where it can improve the performance of image-based search engines and enhance the overall quality of large-scale image datasets handling.

\section*{Acknowledgements}
First author acknowledges the financial support of the Foundation for Science and Technology (FCT, Portugal) in the framework of the grant 2021.07608.BD. Also, both authors acknowledge the financial support of the FCT in the framework of ARISE (DOI 10.54499/LA/P/0112/2020) and R\&D Unit SYSTEC (base UIDB/00147/2020 and programmatic UIDP/00147/2020 funds).

\section*{Author contributions statement}

Investigation and writing the original version of the manuscript: A.B project administration, funding acquisition, supervision and editing the manuscript: R.C. Both authors reviewed the manuscript. 
\section*{Competing interests}
The authors declare no competing interests.
\section*{Data availabilty}
This manuscript does not report data generation or analysis.
\bibliography{sample}

\end{document}